\documentclass[12pt]{article}
\usepackage[utf8]{inputenc}

\usepackage{amsmath}
\usepackage{amscd}
\usepackage{amssymb}
\usepackage{enumerate}
\usepackage{indentfirst}
\usepackage{latexsym}
\usepackage{multicol}
\usepackage{verbatim}
\usepackage{color}

\newcommand{\zl}{\bar\circ}

\def\qed{$\Box$}

\usepackage{theorem}
\numberwithin{equation}{section}

\def\sgn{\, \hbox{\rm sgn}\,  }

\newcommand{\R}{{\mathbb R}}
\newcommand{\bR}{{\mathbb R}}

\def\d{\partial}

\def\ep{\varepsilon}
\def\eps{\eps}
\def\Om{\Omega}

\def\cA{{\mathcal A}}

\newtheorem{defi}{Definition}

\newtheorem{theo}{Theorem}

\theoremstyle{plain}

\newtheorem{proposition}{Proposition}[section]
\newtheorem{lemma}{Lemma}[section]

\theoremstyle{definition}

\begin{document}

\title{\bf Well-posedness of  sudden directional diffusion equations}

\author{Piotr Bogus\l aw Mucha \& Piotr Rybka}
\maketitle

\begin{center}

{Institute of Applied Mathematics and Mechanics \\
University of Warsaw}

{e-mail: p.mucha@mimuw.edu.pl, p.rybka@mimuw.edu.pl}

\end{center}

\date

{\bf Abstract.} Our goal is to establish existence with suitable initial data of solutions to general 
parabolic equation in one dimension, $u_t = L(u_x)_x$, where $L$ is merely a monotone function. We also 
expose the basic properties of solutions, concentrating on  maximal possible
regularity. 
Analysis of solutions with convex initial data explains why we may call them {\it almost classical}.
Some qualitative aspects of solutions, like facets -- flat regions of solutions, are studied too.

\smallskip 

{\it key words:} singular/degenerate parabolic equations, regularity of solutions, almost classical solutions, solution concepts.


\bigskip
\section{Introduction,  problem,  results}

Anisotropic phenomena are plentiful in the natural sciences and technology. The list is long, we name just a few: crystal growth, phase 
transitions, image analysis, models of segmentation.  We will concentrate on mathematical models involving operators like $u_t - \hbox{div}\, L(\nabla u)$, where $L$ 
is a monotone map. Despite many year of research this area is full of challenging problems, the list of contributions is by no means exhaustive,
   \cite{alt-luck},  \cite{mazon2004},  \cite{gurtin},  \cite{BCCN}, \cite{becachan}, \cite{BNP1}, \cite{BNP2}, \cite{BNP3},    \cite{choksi}, \cite{bene-show}, \cite{giacomelli}, \cite{meirmanov}, \cite{taylor}.

In order to set a feasible goal we will concentrate on one dimensional version of these problems, thus we will study:
\begin{equation}\label{i1}
\begin{array}{ll}
 u_t-\frac{d}{dx}L(u_x)=0 \quad\mbox{in } (0,1)  \times (0,T) \equiv
I \times (0,T) =: I_T,\\
u(0,t)=A, \quad u(1,t)=B,
\end{array}
\end{equation}
with a suitable initial datum $u|_{t=0}=u_0$. In general, the boundary data are time dependent.
Our basic assumption is 
$$
L: \bR \to \bR \mbox{ is merely monotone increasing.}
$$
This is a minimal assumption to render (\ref{i1}) formally a parabolic problem.
It also encompasses a rather broad range of phenomena.
 It is enough to say that, we would like to address simultaneously the following equations,
\begin{equation}\label{rn1}
u_t-u_{xx}=0 ,
\end{equation}
\begin{equation}\label{rn2}
 u_t -\frac{d}{dx}(|u_x|+u_x)=0,
\end{equation}
\begin{equation} \label{rn3}
u_t -\frac{d}{dx}(\hbox{sgn}\, (u_x)) =0.
\end{equation}
We see that mere monotonicity of $L$ means that we admit jumps and flat regions of $L$, 
which  represent singular and degenerate parts of the operator.

We notice that in equation (\ref{rn2}) the diffusion turns on only for $u_x>0$, i.e. for areas where $u$ is increasing. There is no evolution  for decreasing parts. This seems to be the simplest model with anisotropic character.
 On the other hand,  in equation (\ref{rn3}) the diffusion is present only in the zone, where  $u_x =0$. 
We know (cf. \cite{brezis}) that $L$ should be understood as a maximal monotone multivalued operator, e.g. $L(p)= \sgn(p)$ is given by
\begin{equation*}
\sgn(p) =\left\{
\begin{array}{ll}
1 & p>0,\\[2pt]
[-1,1] & p=0,\\
-1 & p<0.
\end{array}
\right. 
\end{equation*}
We want to establish in this paper a general existence result and at the same time the maximal possible regularity permitted by the
 generality of our assumptions on $L$.

In principle, one can claim that existence and uniqueness  is well-known. Our
maximal monotone operator $L$ is a subdifferential of a convex function $W$, $L=
\partial W$. Next, we set
\begin{equation}\label{jot}
 J(u) = \left\{
\begin{array}{ll}
\int_I W(u_x)\,dx &\hbox{if the integral is finite},\\
+\infty & \hbox{elsewhere.}
\end{array}
\right.
\end{equation}
Thus, equation (\ref{i1}) is a gradient flow of $J$ in $L_p$ and we can apply
the abstract
nonlinear semigroup theory, see \cite{barbu}, \cite{brezis}. This approach has been implemented by many authors, including \cite{alt-luck}, \cite{mazon2004}, \cite{fukui}. We can call the
solution, we obtain in this way, {\it variational}, because constructing the
Yosida approximation requires solving a variational problem.

However, discussing regularity with this tool
 is cumbersome, because it requires determining the domain of the operator 
$$
\cA:D(\cA) \subset L_2(I) \to L_2(I),
$$
given by
$$
\cA(u) = \frac{\partial}{\partial x}L(u_x).
$$
Of course, $\cA$ has to be properly understood.

An additional complication comes from the fact that the boundary data may be time dependent. 
In any case, there is no common description of  $D(\cA)$. It suffices to consider in (\ref{i1}) $A=0=B,$  and to mention that
if $L(p)=p$, (i.e. we have (\ref{rn1}) in mind),  then $D(\cA) =H^2(I)\cap H^1_0$. 
On the other hand, when $L(p)=\hbox{sgn}\,(p)$,
 then $D(\cA)  {\varsubsetneq}\{u\in AC:\ u|_{\partial I} =0, \ u_x \in BV\}$. In addition, we know that $D(\cA)$  does  not contain all smooth functions,
 e.g. we know that
$x^2\not\in D(\cA)$ (cf. \cite{non}). As a matter of fact a complete description of $D(\cA)$ for that $L$ 
is not available yet. However, it is available  for the multidimensional version of (\ref{rn3}), see e.g. \cite{mazon2010}. It is based on
Anzellotti's formula for integration by parts,
\cite{anzellotti}. However, a direct characterization  of this set for
the one-dimensional problem seems to be missing, a partial result in this direction, for Neumann boundary data, is in \cite{macz}.

Another aspect of the same problem is shown in the papers by Meyer (see \cite{meyer}) and Mucha  (see \cite{mucha}).
These authors showed that for $\cA$ the total variation  operator, i.e. $\cA u= \hbox{div}\, \frac{\nabla u}{| \nabla u|}$ and $f\not\equiv 0$ but small, the only solution to  $\cA u=f$ with the homogeneous Dirichlet boundary condition is $u\equiv 0$.

Therefore, it is not surprising that
the complete discussion of regular of solutions is beyond the
scope of this paper. 

It is worth noticing a development of the viscosity theory for equations like
$$
u_t = a(u_x) \left( W_p(u_x)_x +\sigma\right),
$$
where $W$ is merely convex. This line of research was initiated by Giga and Giga, see \cite{gg-visc}. The point is, this theory is capable of handling faceted solutions too, but facets have to be of positive length. In \cite{miyory} a
comparison principle was proved for such solutions and a number of examples was
presented. However, there is no general existence theory, yet.

Before we discuss the regularity of solutions, we state the basic existence result.

\begin{theo}\label{tw1}
Let us assume that $L(\cdot)$  is monotone increasing, $A,B\in W^1_2(0,T) $,
$u_0 \in L_1(0,1)$, such that $u_{0,x} \in BV(0,1)$, then there 
exists a unique solution $u$ to problem (\ref{i1}), i.e.
$$
u_x \in L_\infty(0,T;BV[0,1]), \qquad u(0,t)=A \mbox{ and } u(1,t)=B
$$
and there exists $\Omega\in L_2((0,1)\times(0,T)$ such that following integral identity is valid
\begin{equation}\label{i2}
(u_t,\phi)_{L_2(0,1)}+(\Omega,\phi_x)_{L_2(0,1)} =0\qquad\hbox{for a.e. }t>0
\end{equation}
and for all $\phi \in C^\infty_0((0,1)\times (0,T))$ and 
\begin{equation}\label{i2-p}
\Omega(x,t) \in L\circ u_x(x,t),
\end{equation}
 where the last term is treated as a composition of multivalued functions.
\end{theo}


In fact, we shall see that the space derivative $\d_x u(\cdot, t) \in BV[0,1]$ for $t>0$ is the
highest possible regularity, 
permitted in such generality. We are able to show that only because we consider a one dimensional problem. Of course, such a general result 
need not be optimal for all choices of $L$.

In order to describe the regularity of solutions we will use here the notion introduced in \cite{non}, \cite{cba}, and developed in \cite{tvn}. Namely, we introduced the notion  
of almost classical solutions. Roughly speaking, $u$, a suitably defined weak solution is {\it almost classical} if it satisfies  the equation
pointwise, except a very small set, see Definition \ref{def:3}. In \cite{non} we
showed that the equation, we studied there, is satisfied except finitely many points. We
will illustrate this notion assuming the initial conditions are convex. This
restriction is for the 
sake of simplicity only, originally in \cite{non}, \cite{tvn} we do not impose
such restrictions. 

The difficulty, related with this notion, is that it is not quite pointwise. In reality we have to give a meaning to the composition of $L$ and $u_x$, when both of them 
are discontinuous. 
It turns out that it is convenient to interpret these two functions as multivalued operators. While is it quite easy to interpret $L$ as a maximal monotone operator, 
more work is necessary if $u_x \in BV$. In this case $u$ is Lipschitz continuous and we may use the notion of Clarke
differential of $u$ to define set-valued $u_x$. The details on the Clarke
differential can be found in the Appendix.
 
Of course, composition of multivalued monotone operators in general need not be monotone. However, we solve this issue 
by introducing $L\bar\circ u_x$, which is a special selection of $L \circ u_x$, see Definition \ref{def:2}. The main advantage is that 
indeed   $L\bar\circ u_x$ is a monotone operator, and on facets
it is even single valued.  
The drawback of this definition is that it is still essentially one-dimensional.

The existence result, Theorem \ref{tw1}, involves passing to the limit in
$L^\epsilon (u^\epsilon_x)$, where quantities indexed by $\epsilon$ are related to an approximate system, when both limiting 
functions $L$ and $u_x$ may
possibly be discontinuous. Thus, it is worthwhile to look more closely at the
limit of $L^\epsilon (u^\epsilon_x)$. The question is to which extent we can
call $\Omega$, the limit of $L^\epsilon (u^\epsilon_x)$, the composition
$L\circ u_x$ or rather $\Omega \subset L\circ u_x$. In principle there is a lot
of choice for $\Omega$, but by 
uniqueness theorem the freedom is rather limited. Finally, it turns out that we can give a 
formula for $\Omega$ that based on an explicitly defined composition $L\bar\circ u_x$.

A few aspects of (\ref{i1}) have been well studied. The most prominent example is the case of uniformly parabolic equation, when $L$ is smooth and
  $\frac{dL}{dp}(p)\ge\epsilon>0$, see \cite{LSU}. The degenerate problems have been studied in \cite{alt-luck} and from  a different perspective in \cite{benedetto}. 
The TV-flow is another extreme case. It has been well studied due to its connection to image processing. There, frequently, 
the domain is rectangular and the solution satisfies the homogeneous Neumann boundary conditions, see \cite{mazon2004}, or it is studied in the whole $\bR^n$, \cite{alter}, \cite{be-jde}. 
 However, the case of Dirichlet data has been also studied, see \cite{mazon2010}. The studies motivated by crystal growth are
 by \cite{fukui}, \cite{G-R}, \cite{koba-giga}, the anisotropic version of the total variation flow (\ref{rn3}) was performed in
 \cite{moll}. The TV flow with constraints was analyzed in \cite{giacomelli}.  Finally, the one dimensional version of the total variation flow was studied in \cite{tvn}
 and \cite{figalli}.
Some  basic questions have been studied for a system with
$L(p)=p+\sgn p$ in \cite{MR-JSP}. However that paper showed just the direction of
our research program, not a complete  general theory.
The present paper  provides the missing existence result.

The complete analysis of (\ref{i1}), taking into consideration all possible types of behavior
of $L$ is beyond the scope of the paper. However, we may relatively easily
discuss the case of $u_0$ being a convex function. The result below guarantees
that convexity is preserved.

\begin{theo}\label{th:con} 
Let the assumptions of Theorem \ref{tw1} are fulfilled. In  addition, $A,B$ are
time independent and  $u_0$ is a convex function. Then for each $t \in (0,T)$
the solution $u(\cdot,t)$ is a convex function, too.
\end{theo}

By Theorem \ref{th:con}, we may restrict our attention to convex data $u_0$. Since at each time instance $u_x(\cdot, t)$ is a monotone 
function, we can extended to a maximal monotone graph. We  define such composition $L\bar\circ u_x$ of maximal monotone graphs
 $L$ and $u_x$ so that $L\bar\circ u_x$ is a monotone operator.

Now we want to distinguish two sets for $f$ and $L$ via a definition below.

\begin{defi}\label{def:1}
Let $u$ be a convex function over $I$, then we set
\begin{equation}
 D_f (u_x) =\Big\{ \bigcup_{k} [a_k,b_k]: u_x|_{[a_k,b_k]}=\theta_k \mbox{ for a
constant } \theta_k \mbox{ and  } a_k < b_k \Big\}.
\end{equation}
Let $L$ be an increasing function over $\R$, then we set
\begin{equation}
 S_L=\{\omega : L(\omega)=[a, b] \mbox{ with } a < b\}.
\end{equation}
\end{defi}
The set $D_f(u_x)$ describes the flat parts of $u_x$ and $S_L$ jump points of $L$.

Now, we are prepared to introduce the definition of  our special composition
$\bar \circ$.

\begin{defi}\label{def:2}
 Let $L$ be a multivalued operator given by (\ref{i2}) and $u_x \in BV(I)$ be a
derivative of a convex function. Then, the multivalued operator
\begin{equation}\label{dd2}
 L\bar\circ u_x
\end{equation}
is defined as follows:

1. If $u_x(x_0) \notin S_L$, then
\begin{equation}\label{dd3}
 L\zl u_x(x_0)=L\circ u_x(x_0).
\end{equation}
In particular, if $u_x(x_0)=[p,q]$, then $L\zl u_x(x_0)=[L(p),L(q)]$.

2. If $u_x(x_0)=\theta \in S_L$ and $L(\theta)=[a,b]$, then:

2.1. if $u_x^{-1}(\theta)=[\xi_-,\xi_+]$ with $\xi_-<\xi_+$ and $\xi_- \neq 0$, $\xi_+ \neq 1$, then
\begin{equation}\label{dd4}
 L\zl u_x(s)=x_k(s-\xi_+)+y_k(s-\xi_-) \mbox{ for } s \in (\xi_-,\xi_+)
\end{equation}
with 
\begin{equation}\label{dd5}
 x_k=\frac{a}{\xi_--\xi_+} \mbox{ and } y_k=\frac{b}{\xi_+-\xi_-};
\end{equation}

2.2. if  $u_x^{-1}(\theta)=\{x_0\}$, then
\begin{equation}\label{dd6}
  L\zl u_x(x_0)=[a,b];
\end{equation}

2.3. if $u_x^{-1}(\theta)=[0,\xi_+]$ with $\xi_+ \in [0,1)$ (or $[\xi_-,1]$ with $\xi_-\in (0,1]$), then we put
\begin{equation}\label{dd7}
  L\zl u_x(s)= b \mbox{ ~~ for ~~} s \in [0,\xi_+] \qquad \Big( L\zl u_x(s)=
a \mbox{ ~~ for ~~} s \in [\xi_-,1] \Big);
\end{equation}

2.4. if $u_x^{-1}(\theta)=[0,1]$, then we put 
\begin{equation}\label{dd8}
 L\zl u_x(s)= \beta \mbox{ ~~ for ~~ } s\in [0,1].
\end{equation}

\end{defi}

By definition, $L\zl u_x =:\sigma$ is an element of a composition of two
multifunctions, $L\circ u_x$. In the examples below, we frequently
specify  a selection $\sigma$, which is absolutely continuous. Thus by
the uniqueness theorem, it is sufficient that the equation (\ref{i1})
is satisfied pointwise a.e. Observe that the composition $L\zl u_x$ is maximal
monotone when $u$ is a convex function.

Finally, we define {\it almost classical} solutions.

\begin{defi}\label{def:3}
{Let us suppose that $u$ is a solution given by  Theorem \ref{tw1}.
We say that $u$ is an
{\underline{almost classical}} solution to system (\ref{i1}) iff the following identity
holds}
\begin{equation}\label{dd1}
 \begin{array}{lcr}
 u_t -\frac{d}{dx} L\bar\circ u_x =0  \quad 
& \mbox{ in } & I \times (0,T)) \setminus E,
\end{array}
\end{equation}
where $\mathcal{H}^1(E)$ is finite.
\end{defi}

The examples we constructed so far (see \cite{non}) show that we cannot expect that $E$ be smaller than stated above.

The definition of the almost classical solutions,  introduced in \cite{tvn}, is based on the
composition $\bar \circ$. We refer to this paper for
details. Here, we present the definition which covers the needs of our
analysis, namely  the case of $u$ being a solution with convex initial data. We shall emphasize that the set of irregular points is at most countable for $t>0$,
but if $L$ suffers just a finite number of jumps, or more precisely, jumps are isolated, then the number of irregular points is finite for 
each $t>0$.

The final result clarifies the importance of solutions to (\ref{i1}) with convex data.

\begin{theo}\label{th:4}
Solutions with convex initial data coming from Theorem \ref{tw1} are indeed almost classical.
\end{theo}

The rest of the paper is devoted to the proofs of these results. In the next
section Theorem \ref{tw1} is shown, in Section 3 we
 prove results for convex solutions: Theorems \ref{th:con} and \ref{th:4}.  Some parts of proofs are stated very precisely, although 
they belong to the classical/well known theory. But they are important to explain the need of redefining the meaning of solutions as well as to initiate the revision of regularity of solutions.

\section{Existence and uniqueness of variational solutions}




In this section we prove Theorem \ref{tw1}.
 We proceed in a standard way. The monotone operator $L$, will be substituted 
by a strictly monotone, single valued smooth function
\begin{equation}\label{i4}
L^\epsilon(p)=(L\ast \pi^\epsilon)(p) + \epsilon p,
\end{equation}
where $\pi^\epsilon$ is a smooth non-negative approximation of the Dirac delta.

Now, we construct an approximating systems. For a given $\epsilon>0$ we examine
\begin{equation}\label{i5}
\begin{array}{lcr}
u^\epsilon_t-\frac{d}{dx}L^\epsilon(u^\epsilon_x)=0 & \mbox{ in } & (0,1)\times (0,T),\\[5pt]
u^\epsilon(0,t)=A^\epsilon, \quad u(1,t)=B^\epsilon & \mbox{ for } & t\in (0,T), \\[5pt]
u^\epsilon|_{t=0}=u^\epsilon_0 & \mbox{ on } & (0,1).
\end{array}
\end{equation}
Data  $A^\epsilon,B^\epsilon,u_0^\epsilon$ are so regularized that they are smooth functions on 
their domains and the consistency conditions are satisfied.

Due to the classical theory, see \cite{LSU},  system (\ref{i5}) admits a unique
smooth solution, at least 
locally in time. Hence,  we concentrate our attention on finding suitable a priori bounds, which 
allow us to pass to the limit with $\epsilon \to 0$. Existence and  uniqueness will be established 
on a given interval $[0,T]$. 

Let us formulate the basic energy estimates guarantying the existence solution to the
approximating system, stipulated by Theorem  \ref{tw1}. 
For this purpose we reduce the given problem to one with homogeneous boundary
conditions. We set 
\begin{equation}\label{i6}
u^\epsilon- [(B^\epsilon(t)-A^\epsilon(t))x+A^\epsilon(t)] =: v^\epsilon.
\end{equation}  
Let us subtract 
$[(B^\epsilon(t)-A^\epsilon(t))x+A^\epsilon(t)]_t 
-\frac{\partial}{\partial x}L^\epsilon([(B^\epsilon(t)-A^\epsilon(t))x+A^\epsilon(t)]_x)$ from both 
sides of (\ref{i5}). Then, we  have
\begin{equation}\label{i6a}
\begin{array}{lcr}
v^\epsilon_t-\frac{d}{dx}[L^\epsilon(u^\epsilon_x)-
L^\epsilon(B^\epsilon(t)-A^\epsilon(t))]=&&\\
-[((B^\epsilon)'(t)-(A^\epsilon)'(t))x+(A^\epsilon)'(t)] + 
\frac{d}{dx} L^\epsilon(B^\epsilon(t)-A^\epsilon(t)) & \mbox{ in } & (0,1)\times (0,T),\\[6pt]
v^\epsilon(0,t)=0, \quad u^\epsilon(1,t)=0 & \mbox{ for } & t\in (0,T), \\[6pt]
v^\epsilon|_{t=0}=u^\epsilon_0-[(B^\epsilon(0)-A^\epsilon(0))x+A^\epsilon(0)] & \mbox{ on } & (0,1).
\end{array}
\end{equation}
Next, we test the equation (\ref{i6a}) with $ v^\epsilon$.
\begin{multline}\label{i7}
\frac{d}{dt} \frac 12 \int_I (v^\epsilon)^2 dx 
+ 
\int_I [L^\epsilon(u^\epsilon_x)-L^\epsilon(B^\epsilon(t)-A^\epsilon(t))]
[u^\epsilon_x-(B^\epsilon(t)-A^\epsilon(t))] dx
=\\
\int_I \Big[-[((B^\epsilon)'(t)-(A^\epsilon)'(t))x+(A^\epsilon)'(t)] 
+ \frac{d}{dx} L^\epsilon(B^\epsilon(t)-A^\epsilon(t))\Big]v^\epsilon dx .
\end{multline}
We notice that for fixed, positive $\epsilon$ the value of 
$\frac{d}{dx} L^\epsilon(B^\epsilon(t)-A^\epsilon(t))$
is clearly zero, because $L^\epsilon$ is single valued.

Identity (\ref{i7}), combined with monotonicity of $L^\epsilon$, leads to the following estimate
\begin{equation*}
 \frac{d}{dt} \frac 12 \int_I (v^\epsilon)^2 dx \le 
\int_I [-[((B^\epsilon)'(t)-(A^\epsilon)'(t))x+(A^\epsilon)'(t)] v^\epsilon dx.
\end{equation*}
Taking into account Gronwall inequality yields a desired
$\epsilon$-independent estimate,
\begin{multline}\label{i8}
\sup_{0\leq t \leq T} \|v^\epsilon(t)\|^2_{L_2(I)} 
\leq \left(\|v^\epsilon_0\|_{L_2(I)}
+\int_0^T\|((A^\epsilon)'(t)-(B^\epsilon)'(t))x-(A^\epsilon)'(t)\|^2_{L_2(I))} dt \right)e^T.
\end{multline}

Let us notice that if $A$ and $B$ are time independent, then we will obtain a better 
estimate in place of (\ref{i8}), i.e.
\begin{equation}\label{i8i}
 \sup_{0\leq t \leq T} \|v^\epsilon\|^2_{L_2(I)} 
\leq \|v^\epsilon_0\|_{L_2(I)}.
\end{equation}

Now, we would like to find a bound on the second space derivative of $u^\epsilon$. After having
differentiated (\ref{i5}) twice with respect to $x$, we will have,
\begin{equation}\label{i9}
\begin{array}{lcr}
u^\epsilon_{xxt}-\frac{d^2}{dx^2}(\frac{d}{dp} L^\epsilon(u^\epsilon_x)u_{xx})=0 & 
\mbox{ in } & I_T,\\[6pt]
\frac{d}{dp} L^\epsilon(u^\epsilon_x)u_{xx}(0,t)=A^\epsilon_t, \quad 
\frac{d}{dp} L^\epsilon(u^\epsilon_x)u_{xx}(1,t)=B^\epsilon_t \qquad & \mbox{ for } & t\in (0,T), \\[6pt]
u^\epsilon_{xx}|_{t=0}=u^\epsilon_{0,xx} & \mbox{ on } & (0,1).
\end{array}
\end{equation}
The boundary data are obtained directly from the system (\ref{i5}) by the analysis of the boundary 
values of $u^\epsilon_t$. 
By the analysis of integrals over set $\{u^\epsilon_{xx}>0\}$ and  $\{u^\epsilon_{xx}<0\}$ we shall 
momentarily show a form of the maximum principle
\begin{equation}\label{i10}
\|u^\epsilon_{xx}\|_{L_\infty(0,T;L_1(0,1))} \leq 
\|u^\epsilon_{0,xx}\|_{L_1(0,1)}.
\end{equation}

Let us fix $\epsilon$. Despite smoothness of  $u^\epsilon$ we cannot claim that the 
set $\{ u^\epsilon_{xx}=0\}$ 
is regular. However, Sard theorem gives us a sequence 
$\sigma_k \to 0$ such that, indeed, each $\{u^\epsilon_{xx} = \sigma_k\}$ is
regular.

Let us fix $\sigma_k$. In order to  examine  set $\{u^\epsilon_{xx}(x) > \sigma_k\}$ 
we look at (\ref{i9}) written in the following form,
\begin{equation}\label{i11}
\begin{array}{lcr}
(u^\epsilon_{xx}-\sigma_k)_t
-\frac{d^2}{dx^2}(\frac{d}{dp} L^\epsilon(u^\epsilon_x)(u_{xx}-\sigma_k))=
\sigma_k\frac{d^2}{dx^2}\frac{d}{dp} L^\epsilon(u_x^\epsilon) 
& \mbox{ in } & (0,1)\times (0,T),\\[6pt]
\frac{d}{dp} L^\epsilon(u^\epsilon_x)u_{xx}(0,t)=A^\epsilon_t, \quad 
\frac{d}{dp} L^\epsilon(u^\epsilon_x)u_{xx}(1,t)=B^\epsilon_t & \mbox{ for } & t\in (0,T), \\[6pt]
u^\epsilon_{xx}|_{t=0}=u^\epsilon_{0,xx} & \mbox{ on } & (0,1).
\end{array}
\end{equation}
Let us integrate equation (\ref{i11}${}_1$) over $\{u^\epsilon (x)_{xx} > \sigma_k\}$. 
Noticing that 
$$
\frac{d}{dt}\int_I (u^\epsilon (x)_{xx}-\sigma_k)_+  dx= 
\int_{\{u^\epsilon_{xx} > \sigma_k\}}(u^\epsilon (x)_{xx}-\sigma_k)_t dx
$$ 
yields after integration by parts, that
\begin{multline}\label{i12}
\frac{d}{dt} \int_{ \{u^\epsilon_{xx} > \sigma_k\}} (u^\epsilon_{xx}-\sigma_k)dx 
- \int_{\d \{u^\epsilon_{xx}>\sigma_k\}}
n_{\d \{u^\epsilon_{xx}-\sigma_k\}}(x)\frac{d}{dx} 
(\frac{d}{dp} L^\epsilon(u^\epsilon_x)u_{xx}-\sigma_k)d \bar\mu
\\
=\sigma_k \int_{\{u^\epsilon_{xx} > \sigma\}}
\frac{d^2}{dx^2}\frac{d}{dp}  L^\epsilon(u_x^\epsilon)dx.
\end{multline}
Here, we denote the atomic measure on $\d \{u^\epsilon_{xx}>\sigma_k\}$ by $d\bar\mu$, 
the term $n_{\d \{u^\epsilon_{xx}>\sigma_k\}}(x)$ is the normal vector, here it is $\pm 1$, 
indicating the orientation. The boundary integral is well-defined due to regularity of
$\d \{u^\epsilon_{xx}>\sigma_k\}$. Moreover, we can bound it too. Indeed,
if for $x_0 \in \d \{ u^\epsilon_{xx}>\sigma_k\}$ and $n_{\d \{u^\epsilon_{xx}>\sigma_k\} }(x_0)=-1$,
then for $x\in (x_0,x_0+\delta)$ we are ensured that
$\frac{d}{dp} L^\epsilon(u^\epsilon_x)(u^\epsilon_{xx}-\sigma_k)>0$, 
hence $\frac{d}{dx} (\frac{d}{dp}  L^\epsilon(u^\epsilon_x)(u_{xx}-\sigma_k))(x_0)\geq 0$. 
By the same token, if 
$$n_{\d \{u^\epsilon_{xx}>\sigma_k \} }(x_0)=1,$$
then for $x\in (x_0-\delta,x_0)$ we are ensured that
$\d L^\epsilon(u^\epsilon_x)(u_{xx}-\sigma_k)>0$. As a result, 
$\frac{d}{dx}  (\d L^\epsilon(u^\epsilon_x)(u_{xx}-\sigma_k))(x_0)\leq 0$. 
Thus, we conclude that 
\begin{equation}\label{i13}
- \int_{\d \{u^\epsilon_{xx}>\sigma_k\}}
n_{\d \{u^\epsilon_{xx}>\sigma_k\}}(x)
\frac{d}{dx}  (\frac{d}{dp} L^\epsilon(u^\epsilon_x)(u_{xx}-\sigma_k))d\bar\mu\geq 0.
\end{equation}
This is why we arrive at
\begin{equation}\label{i14}
\sup_{0\leq t\leq T} \int_I (u^\epsilon_{xx} -\sigma_k)_+ dx 
\leq \int_I (u^\epsilon_{0,xx} -\sigma_k)_+ dx +
\sigma_k \int_0^T \int_I |\frac{d^2}{dx^2} (\frac{d}{dx} L^\epsilon(u^\epsilon_x)) |dxdt.
\end{equation}
We recall that sequence $\sigma_k$ was chosen for fixed $\epsilon>0$. We are allowed to pass 
to the limit with $\sigma_k$. Since the solutions as smooth, the integral in the right-hand-side (rhs)
is well defined and finite, so after taking the limit $\sigma_k \to 0$ we get
\begin{equation}\label{i15}
\sup_{0\leq t\leq T} \int_I (u_{xx})_+ dx \leq \int_I (u_{0,xx})_+ dx.
\end{equation}
Repeating this argument for $(u^\epsilon_{xx})_-=\max \{-u^\epsilon_{xx},0\}$ yields
\begin{equation}\label{i16}
\sup_{0\leq t\leq T} \int_I |u_{xx}^\epsilon |dx \leq \int_I |u_{0,xx}^\epsilon | dx.
\end{equation}

We aim at using the lower semiconituity of the $BV$ semi-norm with respect to $L_1$ convergence. For this purpose we have to establish that $u^\ep_{x}$ tends to $u_{x}$ as $\ep \to 0$. 
In order to capture this we need a bound  on $u^\epsilon_t$. We go back
to the equation and together with the information implied by (\ref{i5}) we obtain
\begin{equation}\label{i17}
\|u_t^\epsilon\|_{L_\infty(0,T;W^{-1}_\infty(I))} \leq 
\|L^\epsilon(u^\epsilon_x)\|_{L_\infty(0,T;L_\infty(I))}\le M.
\end{equation}
The rhs of (\ref{i17}) is uniformly estimated, because by (\ref{i16})
$u^\epsilon_x$ is pointwisely bounded and from the definition of $L^\epsilon$
we deduce that $L^\epsilon( [-1-\|u^\epsilon_x\|_{L_\infty}, 1+\|u^\epsilon_x\|_{L_\infty}])$ 
is a bounded set too.

On the way, we would like to obtain some information in the H\"older class, which is of independent interest. For this purpose, we have to modify a bit the sequence. Take the extension operator
onto the whole line $\R$, such that the function will be compactly supported and the considered regularity will be preserved too.
Namely, we take $U^\epsilon = Eu^\epsilon$. Next, we define
\begin{equation}\label{i18}
V^\epsilon=(1-\d^2_x)^{-1/2}U^\epsilon.
\end{equation}
Estimates (\ref{i16}) and (\ref{i17}) combined yield the following bound for any finite $p>1$,
\begin{equation}\label{i19}
V^\epsilon_t \in L_p(0,T;L_p(\R)),\quad V^\epsilon_{xx} \in L_p(0,T;L_p(\R)),
\qquad\mbox{ i.e. } 
V^\epsilon \in W^{2,1}_p(\R \times (0,T)).
\end{equation}
The embedding yields
$W^{2,1}_p\subset C^{2\alpha,\alpha}(\R \times (0,T)) 
\mbox{~~with~~} \alpha = 1-\frac{3}{2p}>\frac 12 $, provided that for $p>3$,
\cite[Chap. XVIII]{BIN}.
%
 Going back to $u^\epsilon$ we get
 \begin{equation}\label{i22}
u^\epsilon \in C^{2\alpha -1,\alpha -\frac{1}{2}}.
\end{equation}
This clarifies the strong/pointwise convergence of $u$ itself, 
but also we need  information about the strong convergence  of $u^\epsilon_x$.
From (\ref{i16}) we deduce an extra regularity due to 
$u^\epsilon_x \in L_p(0,T;W^{1/p}_p(I)).$
We observe that $L_\infty(0,T;W^2_1(\R)) \subset L_p(0,T;W^{1+1/p}_p(\R))$,
because $\frac 12 (1-\frac 1p) + \frac{1+1/p}{2}=1$, \cite[Chap.XVIII]{BIN}.
Thus, due to the definition of $V^\epsilon$ we obtain that
\begin{equation}\label{i24}
V^\epsilon \in L_p(0,T; W^{2+1/p,1}_p(\R)) 
\mbox{  and    } 
V^\epsilon_{xxt} \in L_p(0,T; W^{-2}_p(\R))\mbox{ ~~  for } 1< p<\infty.
\end{equation}
We can restate (\ref{i24}) in terms of $U$ to obtain
\begin{equation}
 U^\epsilon_x \in L_p(0,T;W^{1/p}_p(\R)) \mbox{ ~~ and ~~ } U^\epsilon_{xt} \in L_p(0,T;W^{-1/p}_p(\R)).
\end{equation}
Since  $W^{1/p}_p \subset L_p \subset W^{-1}_p$ and
%
%
the Aubin-Lion theorem applied to $u^\ep$ defined over $I$ yields (up to a subsequence $\epsilon_k \to0$)
%
\begin{equation}\label{i26}
u^\epsilon_{x} \to u^*_{x} \mbox{ strongly in } L_p(0,T;L_{p}(I)).
\end{equation}
%
Passing to the limit $\epsilon \to 0$, while using (\ref{i16}) and the lower semiconituity of the $BV$ seminorm yields, see \cite{ziemer}, 
\begin{equation}\label{i27}
{\rm ess\,} \sup_{0\leq t\leq T} \|u_{x}(\cdot,t)\|_{BV(I)} \leq \| u_{0,x}\|_{BV(I)}.
\end{equation}
Due to (\ref{i22}) we control the solution at the boundary and 
\begin{equation*}
 u^\epsilon|_{\d I} \to u|_{\d I}  \mbox{ strongly in } C[0,T],
\end{equation*}
and of course the boundary conditions $(\ref{i1})_2$ are fulfilled.

Once we established (\ref{i27}) we want to use it to study
$u^\epsilon_t$ and the composition $L^\epsilon(u^\epsilon_x)$. Our goal is to
show that for
almost all $t$, $L^\epsilon(u^\epsilon_x(\cdot,t))$ has a pointwise limit denoted by $\Om(\cdot,t)$. 
The problem is lack of continuity of $L$.

Let us first set
$g^\epsilon(x,t) = (B^\epsilon (t) - A^\epsilon(t)) x + A^\epsilon(t),$
then $v^\epsilon$ satisfies
\begin{equation}\label{r1}
 v^\epsilon_t - \frac{\d}{\d x} L^\epsilon(u^\epsilon_x) = -g^\epsilon_t.
\end{equation}
Now, we multiply (\ref{r1}) by $v^\epsilon_t$ and integrate over $\Om_T$. After integration by
parts the second term on the l.h.s. and using that $v^\epsilon|_{\d I} =0$ we obtain,
\begin{equation}\label{r2}
 \int_I\int_0^T [(v^\epsilon_t)^2 + L^\epsilon(u^\epsilon_x) v^\epsilon_{xt}] dx dt =
- \int_I\int_0^T g^\epsilon_t v^\epsilon_{t} dx dt.
\end{equation}

Let us now introduce $W^\epsilon$ as a primitive function of $L^\epsilon$,
$$
W^\epsilon (p) = \int_0^p L^\epsilon(r)\, dr +C,
$$ 
where constant $C$ is so chosen that for all 
$p\in [-1-\|u^\epsilon_x\|_{L_\infty}, 1+\|u^\epsilon_x\|_{L_\infty}]$ we have
$W(p)\ge 0$. Monotonicity of $L^\epsilon$ implies convexity of $W^\epsilon$.
Hence, due
to the definition of $v^\epsilon$ we obtain
\begin{multline*}
 \frac{1}{2} \int_I\int_0^T (v^\epsilon_t)^2 dxdt+
\int_I W^\epsilon(u_x) dxdt \\
\le \int_I W^\epsilon(u_{0,x}) dxdt+ 
\frac{1}{2} \int_I\int_0^T (g^\epsilon_t)^2 dxdt  +\int_I\int_0^T
L^\epsilon(u^\epsilon_x)
g^\epsilon_t dxdt.
\end{multline*}
The definition of $v^\epsilon$ yields,
\begin{equation*}
 \int_I\int_0^T (u^\epsilon_t)^2 dxdt
\le \int_I\int_0^T (g^\epsilon_t)^2 dxdt + 4\int_I W^\epsilon(u_{0,x}) dxdt +
\int_I\int_0^T (L^\epsilon(u^\epsilon_x))^2 dxdt.
\end{equation*}
Boundedness of $\| u^\epsilon_x \|_{L_\infty}$ implied by (\ref{i27}) yields that the last term on
the r.h.s. above is finite. Thus, we have shown that
\begin{equation}\label{r3}
 \| u^\epsilon_t,u^\epsilon_x \|_{L_2(\Omega_T)} \le
C(\|A^\epsilon\|_{W^1_2(0,T)}, \|B^\epsilon\|_{W^1_2(0,T)}, L, \|u_{0,x}\|_{BV}),
\end{equation}
where $C$ is independent of $\epsilon$. We  now conclude that
$$
u^\epsilon \rightharpoonup u\qquad\hbox{in } W^{1,1}_2(I_T).
$$
At the same time, by (\ref{i17}) $L^\epsilon(u^\epsilon_x)$ is bounded in $L_2(I_T)$, hence after 
selecting a subsequence
\begin{equation}\label{r3-5}
L^\epsilon(u^\epsilon_x) \rightharpoonup \Omega\qquad\hbox{in } L_2(I_T).
\end{equation}
We will show that the convergence is better.

\begin{lemma}
 The limit functions $u$ and $\Omega$ fulfill (\ref{i2}), i.e. for almost all $t>0$ the identity $u_t =\Omega_x$ holds in the sense of distributions on $(0,1)$.
\end{lemma}

{\bf Proof.} We notice that for a fixed $\varphi\in C^\infty_0(I)$ and any
$t_0\in(0,T)$, $\tau>0$ expression
$$
\int_{t_0-\tau}^{t_0+\tau}\int_I u^\epsilon_t \varphi dxdt
$$
is a continuous functional over $W^1_2(I_T)$. Then, weak
convergence of $u^\epsilon_t$ implies that 
$$
\int_{t_0-\tau}^{t_0+\tau}\int_I u^\epsilon_t \varphi dxdt
\to
\int_{t_0-\tau}^{t_0+\tau}\int_I u_t \varphi dxdt.
$$
Using the equation (\ref{i5}) and (\ref{r3-5}) we shall see that
$$
\int_{t_0-\tau}^{t_0+\tau}\int_I u^\epsilon_t \varphi
=\int_{t_0-\tau}^{t_0+\tau}\int_I
 (L^\epsilon (u^\epsilon_x) )_x\varphi=
-\int_{t_0-\tau}^{t_0+\tau}\int_I 
L^\epsilon (u^\epsilon_x) \varphi_x \to \int_{t_0-\tau}^{t_0+\tau}\int_I \Omega \varphi_x.
$$
Thus, for all $t\in(0,T)$, $\tau>0$ and $\varphi\in C^\infty_0(I)$ we have shown that
\begin{equation}\label{rn5}
 -\int_{t_0-\tau}^{t_0+\tau}\int_I \Omega \varphi_x = 
\int_{t_0-\tau}^{t_0+\tau}\int_I u_t \varphi.
\end{equation}

Since  $\varphi$ is fixed, then Lebesgue differentiation theorem yields
\begin{equation}\label{rn6}
 -\int_I \Omega(x,t) \varphi_x(x)\,dx = \int_I  u_t(x,t) \varphi(x)\,dx,\quad \hbox{for } a.e. 
\ t\in (0,T).
\end{equation}
The choice of $t$ depends on $\varphi$. We  repeat the argument from 
\cite[Theorem 2.1, page 2292]{non} to claim that 
(\ref{rn6}) holds independently of $\varphi$ on a set of full measure $\mathcal{G}\subset(0,T)$.

Identity (\ref{rn6})  means that $\Omega(\cdot,t)$ has a weak derivative and
$\Omega_x(\cdot,t) = u_t(\cdot,t)$
for almost all $t$. In addition,
$$
\|\Omega\|_{L_2(0,T;H^1(I))} \le \|u_t\|_{L_2(I_T)}.
$$

Using the argument as in  \cite[Proposition 3.2]{non} we come to the
following conclusion.
\begin{lemma}\label{lr0}
The limiting function $\Omega$ is a section of the subdifferential $\partial J$ at $u$ (functional $J$ is defined in (\ref{jot})), i.e.
\begin{equation}\label{r4}
\int_I W(u_x+h_x) - W(u_x) dx \ge \int_I \Omega h_x dx.
\end{equation}
\end{lemma}
{\it Proof.} Because of the convexity of functions $W^\epsilon$,
we have the inequalities,
$$
\int_I (W^\epsilon (u^\epsilon_x + h_x) - W^\epsilon (u^\epsilon_x))\,dx \ge 
\int_I \frac{d}{dp}W^\epsilon (u^\epsilon_x))h_x\,dx ,
$$
where $h\in C^\infty_0(I)$. We  integrate both sides with respect to $t$, thus we obtain
\begin{equation}\label{rn7}
\frac{1}{2\tau} \int_{t-\tau}^{t+\tau}\int_I 
(W^\epsilon (u^\epsilon_x + h_x) - W^\epsilon (u^\epsilon_x))\,dxdt \ge 
\frac{1}{2\tau} \int_{t-\tau}^{t+\tau}\int_I \frac{d}{dp}W^\epsilon
(u^\epsilon_x))h_x\,dxdt ,
\end{equation}
We  pass to the limit on the r.h.s. of (\ref{rn7}). If we take into account (\ref{i26}), then we can
apply Krasnosielskii theorem, see \cite{krasno}, to the l.h.s. of (\ref{i26}) while keeping in mind
uniform boundedness of 
$\|u^\epsilon_x\|_{L_\infty}$ given by  (\ref{i27}). This results in the following inequality,
$$
\frac{1}{2\tau} \int_{t-\tau}^{t+\tau}\int_I 
(W (u_x + h_x) - W (u_x))\,dx dt \ge 
\frac{1}{2\tau} \int_{t-\tau}^{t+\tau}\int_I \Omega h_x\,dx dt
$$
Using the argument as before we claim that (\ref{r4}) holds for a.e. $t\in(0,T)$.

We  show more,

\begin{lemma}\label{pr1}
{\sl If (\ref{r4}) holds for a fixed $t$ and  all $h\in H^1_0(I)$, then 
$\Omega(x,t)\in \d W \circ u_x(x,t)$ for all $x$, except at most countably many. In other words $\Omega(x,t) \in \circ u_x(x,t)$ a.e.}
\end{lemma}
{\it Proof.} Since $u_x(\cdot,t) \in BV(I)$, then
$u_x$ is continuous except for at most countably many points. Let us suppose
that $u_x(\cdot,t_0)$ is continuous at $(x_0,t_0)$ and $p\in \bR$ is arbitrary.
We take $h^\epsilon \in C^\infty_0(I)$ such that
$$
\hbox{supp}\, h^\epsilon\subset (x_0-\epsilon, x_0+\epsilon),\quad
h^\epsilon_x(x_0) =p, \quad | h^\epsilon_x | \le M,
$$
where $M$ is independent of $\epsilon$.

Now, we notice that for a function $g$ continuous at $x_0$ we have
$\frac{1}{2\epsilon} \int_{x_0-\epsilon}^{x_0+\epsilon} g(s)\,ds \to g(x_0).$
This observation and (\ref{r4}) combined yield
$$
W(u_x(x_0,t_0)+p) - W(u_x(x_0,t_0) )\ge \Omega(x_0,t_0) p,\quad\hbox{ a.e. }x_0, t_0,
$$
as desired. \qed

%
%
%

\begin{proposition}\label{pr-u}
Variational solutions to (\ref{i2}) are unique.
\end{proposition}

{\it Proof.}
In order to show the uniqueness of our solutions it is enough to work in the $L_2$-framework. We first notice that the difference of two solutions is an admissible test function. This is so, because the difference vanishes on the parabolic boundary of $I_T$.
Take two solutions to (\ref{i2}),  say $u_1$ and $u_2$, then we have
\begin{equation}\label{u1}
(u_{1,t}-u_{2,t}, \phi) + (\Omega_1-\Omega_2,\phi_x)=0 \mbox{ in } I \times (0,T)
\end{equation}
  and 
$u_1(0,t)-u_2(0,t)=0$  and $u_1(1,t)-u_2(1,t)=0$ for $t\in (0,T)$.
Thus, we insert it into (\ref{u1}), 
\begin{equation}\label{u3}
\frac{d}{dt} \frac 12 (u_1-u_2,u_1-u_2) + (\Omega_1-\Omega_2,u_{1,x}-u_{2,x})=0.
\end{equation}
By   Lemma \ref{pr1}, we deduce that 
$(\Omega_1-\Omega_2,u_{1,x}-u_{2,x}) \geq 0$,
since 
the initial data are the same, so (\ref{u3}) yields $\frac{d}{dt}\| u_1-u_2\|^2 \le 0$, i.e.
$u_1 \equiv u_2$ on the full interval of the existence. \qed

\section{Convexity of solutions}

In the previous section we showed existence and uniqueness of solutions to (\ref{i1}). 
Here, we prove results for convex solutions, i.e. Theorems \ref{th:con} and \ref{th:4}.
We proceed as in the proof of (\ref{i10}).

\smallskip 

 {\it Proof of Theorem \ref{th:con}. } 
We aim at proving that for $t>0$,  solution $u(\cdot,t)$ is a convex function.
Since $u$ is obtained as a limit of approximating solutions $u^\epsilon$ to
(\ref{i5}), it is  sufficient to show that for all $t>0$ and $\epsilon >0$
functions $u^\epsilon(\cdot, t)$ are convex.

Let us consider system (\ref{i5}).
Since $u^\epsilon$ at the boundary is fixed and we assumed that $A,B$ are time
independent, then $u^\epsilon_t=0$ at
$\d I$, so $\frac{d}{dp}L^\epsilon(u^\epsilon_x) u^\epsilon_{xx}=0$ at $\d I$, too. As a result, we find that $u^\epsilon_{xx}=0$ at $\d I$.
Making  use of the fact that $u^\epsilon$ are smooth, we are allowed
to differentiate twice eq. (\ref{i5}). This yields
\begin{equation}\label{c1}
 \begin{array}{lcr}
u^\epsilon_{xxt}-\frac{\partial^2}{\partial x^2}[\frac{d}{dp}L^\epsilon(u^\epsilon_x) u^\epsilon_{xx}]=0 & \mbox{ in } & I \times  (0,T),\\[4pt]
u^\epsilon_{xx}=0 & \mbox{ at } & \d I \times (0,T), \\[4pt]
u^\epsilon_{xx}|_{t=0}=u^\epsilon_{0,xx}\geq 0 & \mbox{ at} & I.
\end{array}
\end{equation}

Although $u^\epsilon$ is smooth, we cannot claim that the sets $\{(x,t)
\in I\times [0,T] : u^\epsilon_{xx}(x,t)=0\}$ 
are regular. However, the Sard theorem guarantees us that there exists
a sequence $\delta_k \to 0^-$,  such that sets  
$\{(x,t) \in I\times [0,T] : u^\epsilon_{xx}(x,t)=\delta_k\}$ are
regular submanifolds, where $\delta_k <0$. 
We follow the notation $f_-=\max\{-f,0\}$.
Next, we restate (\ref{c1}) as follows
\begin{equation}\label{c22}
 (u^\epsilon_{xx}-\delta_k)_t-\frac{\partial^2}{\partial
    x^2}[\frac{d}{dp}L^\epsilon(u^\epsilon_x) (u^\epsilon_{xx}-\delta_k)]=
\delta_k \frac{\partial^2}{\partial x^2}\frac{d}{dp}L^\epsilon(u^\epsilon_x).
\end{equation}
Next, we integrate (\ref{c22}) over the set $A_k=\{(x,t) \in I\times [0,T] : u^\epsilon_{xx}(x,t)<\delta_k\}$
$$
\int_{\d A_k} (u_{xx}^\epsilon -\delta_k) n_t d\sigma - \int_{\d A_k} 
\frac{\partial}{\partial x}[\frac{d}{dp}L^\epsilon(u^\epsilon_x) u^\epsilon_{xx}-\delta_k)] n_x d\sigma = \delta_k \int_{A_k}
\frac{\partial ^2}{\partial x^2}\frac{d}{dp}L^\epsilon(u^\epsilon_x).
$$
Here, $(n_x,n_t)$ is the outer normal vector to $A_k$.

Since the set $\d A_k$ is regular and the function $u_{xx}^\epsilon$ equals $ \delta_k$ on $\d A_k \setminus I\times
\{0,T\}$, then we get
$\int_{\d A_k} (u_{xx}^\epsilon -\delta_k) n_t d\sigma = 0,$
moreover,
$\int_{\d A_k} \frac{d}{dx}[\frac{d}{dp}L^\epsilon(u^\epsilon_x) u^\epsilon_{xx}-\delta_k)] n_x d\sigma <0,$
because of $\frac{d}{dp}L^\epsilon(u^\epsilon_x)  >0$ and $\d I \cap (\d A_k \cap I\times \{t\})=\emptyset$ for $t\in (0,T)$.
The initial datum $u^\epsilon_{0,xx}$ is non-negative and
$\delta_k <0$.
Thus, we obtain
\begin{equation}
\int_I (u_{xx}^\epsilon -\delta_k)_- dx  \leq |\delta_k | \int_0^T
 \int_I |\frac{\partial}{\partial x}(\frac{d}{dp}L^\epsilon(u^\epsilon_x))| dxdt.
\end{equation}
After letting $k \to \infty$ ( $\delta_k \to 0^-$), we get
\begin{equation}
 \int_I (u_{xx}^\epsilon)_- dx =0. \mbox{ ~~ Thus ~~ } u_{xx}^\epsilon(x,t) \geq 0 \mbox{ for } t >0.
\end{equation}

Next, after passing with $\epsilon$ to $0$, we complete the proof of Theorem \ref{th:con}.

\bigskip 

{\it Proof of Theorem \ref{th:4}.} We proceed in a few steps. We know that $u_t\in L_2(I_T)$ hence, for almost all $t$, with respect to one-dimensional Lebesgue measure, we have $u(\cdot, t)\in L_2(I)$. We have also established that for almost all $t$ the equality
\begin{equation}\label{a1-1}
u_t = \Omega_x
\end{equation}
holds in the $L_2$ sense. In particular $\Om\in H^1(I)$.

Here comes the main observation. 

\begin{lemma}
 Let us suppose that $u_0$ is convex and $u$ is a corresponding weak solution to (\ref{i1}), then for almost all $t$ function $\Omega(\cdot, t)$ satisfies $\Omega(x,t)=L\bar \circ u_x(x,t)$. Hence, $\Omega$ is differentiable except at most countably many points.
\end{lemma}
{\it Proof.}
We have to consider the cases specified in the definition of $\bar\circ$. 
The crucial one is when for a number $c$ we have $L(c)=[a,b]$ with $a<b$. 
From Theorem \ref{th:con}, we know that 
\begin{equation}\label{a1}
\{x\in  I :u_x(x,t)=c, t>0\}=[\xi_-(t),\xi_+(t)]
\end{equation}
where $\xi_-(t)<\xi_+(t)$, since $u(\cdot,t)$ is convex.

The monotonicity of $L(\cdot)$ and $u_x(\cdot,t)$ implies that we are able to find two sequences $\{x_n^-(t)\},\{x_n^+(t)\}$ such that
\begin{equation}\label{a2}
x_n^-(t) \to \xi^-(t)^- \mbox{ and } x_n^+\to \xi^+(t)^+
\end{equation}
as well as 
\begin{equation}\label{a3}
L(u_x(x_n^-(t),t)) \to a \mbox{ and } L(u_x(x_n^+(t),t)) \to b.
\end{equation}
Namely, we require that points $x_n^\pm(t)$ are  regular of function $u_x(\cdot,t)$ and 
$u_x(x_n^\pm(t),t)$ are regular of $L(\cdot)$, too.
Then, we look at 
\begin{equation}\label{a4}
u_t=\frac{\partial}{\partial x}\Omega(x,t) \mbox{ on } [x_n^-(t),x_n^+(t)].
\end{equation}

The integration over the marked set yields
\begin{equation}\label{a5}
\int_{x_n^-(t)}^{x_n^+(t)} u_t dx = L(u_x(x_n^+(t),t)) -L(u_x(x_n^-(t),t)).
\end{equation}
The form of the rhs follows from the choice of suitable regular points.
Since we are allowed to pass to the limit with $n$, we get
\begin{equation}\label{a6}
\int_{\xi_-(t)}^{\xi_+(t)} u_t dx = b-a.
\end{equation}

On the other hand, function $u$ restricted to the set $[\xi_-(t),\xi_+(t)]$ is linear, 
thus it has the following form
\begin{equation}\label{a7}
u|_{[\xi_-(t),\xi_+(t)]}=cx+\alpha (t).
\end{equation}
Note that the facet is prescribed  by the value $c$, so the slope is fixed in time, however 
the shift 
given by $\alpha(t)$ may depend on time. Of course, this holds only on the considered set. 
Hence, (\ref{a7}) follows
\begin{equation}\label{a8}
u_t|_{[\xi_-(t),\xi_+(t)]}=\alpha_t(t).
\end{equation}
Combining this  with (\ref{a6}), we get the following relation
\begin{equation}\label{a9}
\alpha_t(\xi_+(t)-\xi_-(t))=b-a.
\end{equation}
By Definition \ref{def:2}, we observe that 
\begin{equation}\label{a10}
L \bar\circ u_x|_{[\xi_-,\xi_+]}=\frac{b-a}{\xi_+-\xi_-} x + \frac{a\xi_+}{\xi_+-\xi_-} -\frac{b\xi_-}{\xi_+-\xi_-}.
\end{equation}
Hence, comparing (\ref{a1-1}), (\ref{a9}) and (\ref{a10}) we found 
\begin{equation}\label{a11}
u_t = \frac{d}{dx} L\bar \circ u_x \mbox{ at } [\xi_-(t),\xi_+(t)].
\end{equation}
On the remaining parts of the domain we have the regular case. By Lemma \ref{pr1}, we are ensured that $\Omega \in L\circ u_x$, and since by 
(\ref{a3}) $L(u_x(x_n^-(t),t))=\Omega(x_n^-(t),t) \to a$ and $L(u_x(x_n^+(t),t))=\Omega(x_n^+(t),t) \to b$, and
$L\zl u_x (\xi^-(t),t)=a$ and $L\zl u_x(\xi^+(t),t)=b$. So, we get 
\begin{equation}
 \Omega(x,t) = L\zl u_x \mbox{ as } H^1(0,1) \mbox{ functions}.
\end{equation}
We should stress that simplicity of the consideration follows from the convexity of 
function $u(\cdot,t)$. 
The proof is done. \qed

\section*{Appendix}
%

In (\ref{i2-p}) we compose two multivalued operators $L$ with $u_x$. Here we 
explain how we can treat $BV$ functions as multi-valued operators.
This is easy for functions which are derivatives, $u_x\in
BV(I)$. It is very useful in the regularity study of solution to
(\ref{i1}). If $u$ and $u_x$ belong to $ BV(I)$, then $u$ is
Lipschitz continuous. Hence, $\frac{d^+u}{dx}$ and $\frac{d^-u}{dx}$
exist everywhere and they differ on at most countable set. Thus, we may set
\begin{equation}\label{depaw}
\partial_x u(s) = \{\tau u_x^- + (1-\tau) u_x^+:\ \tau\in [0,1]\}.
\end{equation}
Under our assumptions on $u$, the set $\partial_x u(x)$ is the Clarke
differential of $u$ and equality holds in (\ref{depaw}) due to
\cite[Section 2, Ex. 1]{chang}. If $u$ is convex, then $\partial_x u$
is the well-known subdifferential of $u$. As a result, if $u_x\in BV$,
then for each $ x_0\in (0,1)$, we have
$\partial_x u(x_0)=[\lim_{x\to x_0^-} u_x(x),\lim_{x\to x_0^+} u_x(x)]_{or}$,
where $[a,b]_{or}=[a,b]$ for $a\leq b$ and $[a,b]_{or}=[b,a]$ for $b>a$.


\bigskip\noindent
{\bf Acknowledgment.} The research for this paper has been supported by NCN grant Nr 2011/01/B/ST1/01197.

\end{document}